\documentclass[11pt,draft]{article}
\usepackage{mathrsfs}
\usepackage{amsfonts}
\usepackage{mathrsfs}
\usepackage{amsmath,amsfonts,mathrsfs,amssymb,color}
\usepackage{indentfirst}

\numberwithin{equation}{section}

\newtheorem{theo.}{\quad\, Theorem}[section]
\newtheorem{defi.}{\quad\, Definition}[section]
\newtheorem{lemm.}{\quad\, Lemma}[section]
\newtheorem{coro.}{\quad\, Corollary}[section]

\begin{document}
\title {Proof of the Bonheure-Noris-Weth conjecture on \\ oscillatory radial solutions of Neumann problems\\
 }
\author{
Ruyun Ma$^{a,*}$\ \ \ \ \ Tianlan Chen$^{b}$ \ \ \  \ Yanqiong Lu$^c$
\\
{\small $^{a,b,c}$Department of Mathematics, Northwest
Normal University, Lanzhou 730070, P R China}\\
}
\date{} \maketitle
\noindent\footnote[0]{E-mail addresses:mary@nwnu.edu.cn(R.Ma),\,chentianlan511@126.com(T.Chen),\,linmu8610@163.com(Y.Lu)
} \footnote[0] {$^*$Supported by the
NSFC (No.11361054), SRFDP(No.20126203110004) and Gansu Provincial National Science Foundation of China (No.1208RJZA258). }

{\small\bf Abstract.} {\small  Let $B_1$ be the unit ball in $\mathbb{R}^N$ with $N \geq 2$. Let $f\in C^1([0, \infty), \mathbb{R})$, $f(0)=0$,
$f(\beta) = \beta, \ f(s)<s\ \text{for}\ s\in (0,\beta), \ f(s)>s\ \text{for}\ s\in (\beta, \infty)$ and $f'(\beta)>\lambda^{r}_k$.
 D. Bonheure,   B. Noris and T.  Weth [Ann. Inst. H. Poincar\'{e} Anal. Non Lin\'{e}aire 29(4) (2012)] proved the existence of nondecreasing, radial positive solutions of the semilinear Neumann problem
$$
-\Delta u+u=f(u) \ \text{in}\ B_1,\ \ \ \  \partial_\nu u=0 \ \text{on}\ \partial B_1
$$
for $k=2$, and they conjectured that there exists a radial solution with $k$ intersections with $\beta$ provided that $f'(\beta) >\lambda^r_k$ for $k>2$. In this paper, we show that the answer is yes.
}
\vskip 3mm


\vskip 3mm

{\small\bf MR(2000)\ \ \ 34B10, \ 34B18}

\baselineskip 20 pt

\section{Introduction}

Let $B_1$ be the unit ball in $\mathbb{R}^N$ with $N \geq 2$.
 Very recently, D. Bonheure,   B. Noris and T.  Weth [1] proved the existence of nondecreasing, radial positive solutions of the semilinear Neumann problem
$$
\left\{
\aligned
-\Delta u+u=f(u)    \ \ \ \ \ \ \ & \text{in}\ B_1,\\
u>0\qquad\qquad\quad \ \ \ \ \           & \text{in}\ B_1,\\
\partial_\nu u=0 \ \ \ \qquad\qquad \ \ \ &\text{on}\ \partial B_1 \\
\endaligned
\right.
\eqno (1.1)
$$
under the assumptions:

\noindent $(f1)$ $f\in C^1([0, \infty), \mathbb{R})$, $f(0)=0$ and $f$ is nondecreasing;

\noindent  $(f2)$  $f'(0)=\lim_{s\to 0^+}\frac {f(s)}s=0$;

\noindent  $(f3)$ $\underset{s\to+\infty}\liminf \frac{f(s)}s>1$;

\noindent  $(f4)$ there exists $\beta > 0$ such that $f(\beta) = \beta$ and
$$f'(\beta)>\lambda^{r}_2.
\eqno (1.2)
$$
Here $\lambda^{r}_k$ is the $k$-th radial eigenvalue of $-\Delta +I$ in the unit ball with Neumann boundary conditions.

It is easy to see that $u\equiv \beta$ is a constant solution of (1.1), and  there exists nonlinearity $f$ satisfying $(f1)-(f3)$ such that the problem (1.1)
only admits this constant solution, see [1, Proposition 4.1]. For the existence of nonconstant
 radial solutions, they obtained the following result by variational argument.

\vskip 3mm

\noindent{\bf Theorem A.} Assume $(f1)-(f4)$. Then there exists at least one nonconstant increasing radial solution of (1.1).

\vskip 3mm

  They raised the question whether it is possible to
construct radial solutions with a given number of intersections with $\beta$ provided that $f'(\beta)$ is sufficiently large. More precisely, they conjectured that there exists a radial solution with $k$ intersections with $\beta$ provided that $f'(\beta) >\lambda^r_k$.

   The purpose of the present paper is to show that the answer to the above question is yes! The proof is based upon  the unilateral global bifurcation theorem [4, 5, 8]. The condition $f'(0) = 0$  and
   the monotonic condition in $(f1)$ seem unduly restrictive.  We shall make the following
assumptions:

\vskip 2mm

\noindent $(A1)$ $f\in C^1([0, \infty), \mathbb{R})$, $f(0)=0$;

\noindent $(A2)$  $f_{+\infty}:=\underset{s\to+\infty}\lim \frac{f(s)}s <\infty$;

\noindent $(A3)$ there exists $\beta > 0$ such that
$$f(\beta) = \beta, \ \ \ f(s)<s\ \text{for}\ s\in (0,\beta), \ \ \ f(s)>s\ \text{for}\ s\in (\beta, \infty)$$ and
$$f'(\beta)>\lambda^{r}_k, \ \ \ \ \ \text{for some} \ k\geq 2;
$$

\noindent $(A4)$ $[f(s+\beta)-(s+\beta)]s>0, \ s\in (-\beta, 0)\cup (0, \infty)$.

\vskip 3mm

   The main result of this paper is the following

\vskip 3mm

\noindent{\bf Theorem 1.1} Assume (A1)-(A3).
\noindent Then for each $j\in \{2, \cdots, k\}$, (1.1) has two nonconstant radial solutions
$u_j^+$ and $u_j^-$ such that $u_j^+ -\beta$ changes sign exactly $k-j+1$ times in $(0, 1)$ and is positive near $0$, and $u_j^- -\beta$ changes sign exactly $k-j+1$ times in $(0, 1)$ and is negative near $0$.
Moreover, if $(A4)$ holds, then $u_2^+$ is decreasing in $[0,1]$ and $u_2^-$ is increasing in $[0,1]$.

\vskip 3mm

For other results on the existence of radial solutions of nonlinear Neumann problems,  see [2, 10, 14].

\vskip 3mm

The rest of the paper is organized as follows. In Section 2 we study the spectrum structure of the linear Neumann  problem
$$
\left\{
\aligned
-\Delta u(x)&=\mu a(|x|) u(x)   \ \ \ \ \qquad \ \ \   \text{in}\ B_1,\\
\ \partial_\nu u&=0 \ \  \qquad\quad\qquad\qquad \ \ \ \text{on}\ \partial B_1, \\
\endaligned
\right.
$$
where $a\in C[0, 1]$ satisfies $a(r)>0$ for $r\in [0, 1]$. In Section 3, we introduce some
functional setting and state some preliminary bifurcation results on abstract operator equations.  Finally in Section 4 we prove our main results on the existence of nonconstant radial solutions by applying the well-known unilateral bifurcation theorem due to Dancer [4, 5].

\vskip 3mm

\section{Eigenvalues of linear eigenvalue problems}

Let us consider
the linear eigenvalue problem
$$
\left\{
\aligned
-\Delta u(x)&=\mu a(|x|) u(x)   \ \ \ \ \qquad \ \ \   \text{in}\ B_1,\\
\ \partial_\nu u&=0 \ \  \qquad\quad\qquad\qquad \ \ \ \text{on}\ \partial B_1, \\
\endaligned
\right.
\eqno (2.1)
$$
where $a\in C[0, 1]$ satisfies
$$a(r)>0,\ \ \ \ \  \ r\in [0, 1].
\eqno (2.2)$$

\vskip 3mm

\noindent{\bf Theorem 2.1} Assume that (2.2) is fulfilled. Then the radial eigenvalues of (2.1) are as follows:
$$0=\mu_0^r<\mu_1^r<\mu_2^r< \cdots\to \infty.
\eqno (2.3)
$$
Moreover, for each $k\in \mathbb{N}^*:=\{0, 1, 2, \cdots\}$, the radial eigenvalue $\mu_k^r$ is simple, and the
 radial eigenfunction $\psi_k$, being regarded as a function of $r$, possesses exactly $k$ simple zeros in $[0, 1]$,
and $\psi_k$ is radially monotone if and only if $k\in \{0, 1\}$.

\vskip 3mm

It is easy to see that Theorem 2.1 is an immediate consequence of the following results on singular
Sturm-Liouville problems.

\vskip 3mm

\noindent{\bf Theorem 2.2} Assume that (2.2) is fulfilled. Then the eigenvalues of the problem
$$
\left\{
\aligned
&-u''(r)-\frac{N-1}r u'(r)&=\mu a(r) u(r),   \ \ \ \ \qquad \ \ \   r\in (0,1),\\
&u'(0)=0=u'(1)\\
\endaligned
\right.
\eqno (2.4)
$$
 are as follows:
$$0=\mu_0^r<\mu_1^r<\mu_2^r< \cdots\to \infty.
$$
Moreover, for each $k\in \mathbb{N}^*$, $\mu_k^r$ is simple, and the
 eigenfunction $\psi_k$ possesses exactly $k$ simple zeros in $[0, 1]$,
and $\psi_k$ is monotone if and only if $k\in \{0, 1\}$.

\vskip 3mm

To prove Theorem 2.2, we need several basic lemmas.

\vskip 3mm

\noindent{\bf Lemma 2.1}
    Assume that $\tilde f\in C([0, \infty)\times [0, \infty))$ is Lipschitz continuous  in $u$ on $[0, \infty)$. Then for given $\zeta\in (0, \infty)$,
the  initial value problem
$$
\left\{
\aligned
&-u''(r)-\frac{N-1}r u'(r)&=\tilde f(r, u(r)),   \ \ \ \ \qquad \ \ \   r\in (0,\infty),\\
&u'(0)=0,\\
&u(0)=\zeta\\
\endaligned
\right.
\eqno (2.5)
$$
  has a unique solution  $u$ defined on $[0, \infty)$.
  Moreover, all of zeros of $u$ are simple.

 \noindent{\bf Proof.} \ According to [15,  Existence and uniqueness Theorem XIII in \S 6 of Chapter II\,], for given $b>0$, the initial value problem
 $$
\left\{
\aligned
&-u''(r)-\frac{N-1}r u'(r)&=\tilde f(r, u(r)),   \ \ \ \ \qquad \ \ \   r\in (0,b],\\
&u'(0)=0,\\
&u(0)=\zeta\\
\endaligned
\right.
$$
 has exactly one solution $u\in C^2[0, b]$. Notice the equation in $(2.5)$ is non-singular for $r\geq b$, it is evident that $u$ can be extended to $[0, \infty)$. Since $\tilde f\in C([0, \infty)\times [0, \infty))$ is Lipschitz continuous  in $u$ on $[0, \infty)$, the uniqueness part
can be deduced by the same method in the Appendix in [11].

 All of zeros of $u$ are simple since for any zero point $\tau$ of $u$,
  the initial value problem
 $$
\left\{
\aligned
&-u''(r)-\frac{N-1}r u'(r)&=f(r, u(r)),   \ \ \ \ \qquad \ \ \   r\in (0,\infty),\\
&u(\tau)=0=u'(\tau)\\
\endaligned
\right.
$$
has only trivial solution $u\equiv 0$. \hfill{$\Box$}

\vskip 3mm

\noindent{\bf Lemma 2.2}  Assume that $a\in C([0, \infty),(0, \infty))$ and there exist two positive constants $a_1$ and $a_2$, such that
$$a_1\leq a(r)\leq a_2,\ \ \ \ \  \ \ \  r\in [0, \infty).
$$
Let $u$ be a solution of the problem
$$
\left\{
\aligned
&-u''(r)-\frac{N-1}r u'(r)&=\mu a(r) u(r),   \ \ \ \ \qquad \ \ \   r\in (0,\infty),\\
&u'(0)=0,\\
&u(0)=\zeta\\
\endaligned
\right.
\eqno (2.6)
$$
 with $\mu>0$ and $\zeta>0$.  Then $u$ has a sequence of zeros  $\{\tau_n\}\subset (0, \infty)$ with
$$\tau_n\to \infty\ \ \ \ \ \text{as} \ n\to \infty.
\eqno (2.7)
$$

\noindent{\bf Proof.} From [15,  XVIII in \S 27 of Chapter VI], the solution $y$ of the initial value problem
$$(r^{N-1}y')'+r^{N-1}y=0, \ \ \ \ y(0)=1, \ \ \ y'(0)=0
\eqno (2.8)
$$
 oscillates.
Denote the zeros of $y$ by $\xi_0<\xi_1<\xi_2<\cdots$. Then
$$\xi_{n+1}-\xi_{n}\to \pi, \ \ \ \ \ \ \ \text{as}\ n\to\infty.
\eqno (2.9)
$$

Let $\gamma_j=\sqrt{\mu a_j}$ for $j=1,2$. Let
$$u_j(r)=y(\gamma_j r), \ \ \ \ \ j=1,2.
$$
Then $u_j$ is the unique solution of the initial value problem
$$(r^{N-1}u')'+\gamma_j^2 r^{N-1}u=0, \ \ \ \ u(0)=1, \ \ \ u'(0)=0,
\eqno (2.10)
$$
and for $j=1,2$, $u_j$ oscillates and it has a sequence of zeros
$\frac{\xi_0}{\gamma_j}<\frac{\xi_1}{\gamma_j}<\frac{\xi_2}{\gamma_j}<\cdots$.
Combining this  with the Sturm-Picone Theorem ( see [15, \S 27 of Chapter VI]), it deduces that
the solution $u$ of (2.6) oscillates.
\hfill{$\Box$}

\vskip 3mm

\noindent{\bf Lemma 2.3} Assume that $a\in C([0, \infty),(0, \infty))$. Let $u$ be a solution of the
problem
(2.6) with $\mu>0$ and $\zeta>0$. Let $r_1,r_2$ be any two consecutive zeros of $u'$ in $[0,\infty)$ with $r_1<r_2$. Then $u$ has one and only one zero in $(r_1, r_2)$.

\noindent{\bf Proof.} $-u''(r)-\frac{N-1}r u'(r)=\mu a(r) u(r)$ can be rewritten as
$$-(r^{N-1}u'(r))'=\mu a(r)r^{N-1}u(r).
\eqno (2.11)$$
Integrating from $r_1$ to $r$, we get
$$u'(r)=-\mu \int^r_{r_1} \big(\frac tr\big)^{N-1}a(t)u(t)dt,
$$
and accordingly,
$$0=u'(r_2)=-\mu \int^{r_2}_{r_1} \big(\frac{ t}{r_2}\big)^{N-1}a(t)u(t)dt,
$$
which implies that $u$ has at least  one zero in $(r_1,r_2)$.

   Suppose on the contrary that $u$ has two zeros $z_1, z_2\in (r_1,r_2)$. Then
   there exists $z^*\in (z_1, z_2)\subset (r_1,r_2)$, such that $u'(z^*)=0$. However,
   this is a contradiction.
\hfill{$\Box$}

\vskip 3mm

\noindent{\bf Lemma 2.4} Assume that $a\in C([0, \infty),(0, \infty))$. Let $u$ be a solution of (2.6) with $\mu>0$ and $\zeta>0$. Let $\tau_1,\tau_2$ be any two consecutive zeros of $u$ in $(0,\infty)$ with $\tau_1<\tau_2$. Then $u'$ has one and only one zero in $(\tau_1, \tau_2)$.

\noindent{\bf Proof.} Obviously, $u'$ has at least one zero $r_1$ in $(\tau_1, \tau_2)$.

Without loss of generality, we may assume that $u(r)>0$ in $(\tau_1, \tau_2)$.
It follows from
$$-u''(r)-\frac{N-1}r u'(r)=\mu a(r) u(r)
\eqno (2.12)
$$
that
$$u''(r_1)<0,
$$
which implies that $u$ is concave up near $r=r_1$.

Suppose on the contrary that there exists $r_*\in (\tau_1, \tau_2)$ with $r_*\neq r_1$ such that $u'(r_*)=0$.
Then by the same argument, we get
$$u''(r_*)<0
$$
This together with $u''(r_1)<0$ imply that
there exists $\hat r\in (\min\{r_1, r_*\},\max\{r_1, r_*\})$
such that $u$ attains a local minimum at $\hat r$, and
$$u(\hat r)>0, \ \ \ \ u'(\hat r)=0, \ \ \ \ u''(\hat r)\geq 0,
$$
which contradicts (2.12). Therefore, $u'$ has only one zero in $(\tau_1, \tau_2)$.
\hfill{$\Box$}

\vskip 3mm

\noindent{\bf Lemma 2.5}  Assume that $a\in C([0, \infty),(0, \infty))$ with $a(r)\geq a_0>0$
in $(0, \infty)$. Let $u$ be a solution of (2.6) with $\zeta>0$ and $\mu>0$. Let $\tau_k(\mu)$ and $r_k(\mu)$  be the $k$-th positive zero of $u$ and $u'$, respectively.
 Then

 \noindent(1) For given $k\in \mathbb{N}$, $\tau_k(\mu)$ is strictly decreasing in $(0, \infty)$;

 \noindent(2) For given $k\in \mathbb{N}$, $r_k(\mu)$  is strictly decreasing in $(0, \infty)$.

\vskip3mm

\noindent{\bf Proof.} (1) For fixed $k>1$, $\tau_k(\mu)$ is strictly decreasing in $\mu$, which
is an immediate consequence of the well-known Sturm Separation Theorem
[15, P. 272] since the differential equation
$$-u''-\frac{N-1}r u'=\mu a(r) u$$
is non-singular for $r\geq \tau_1(\mu)$. So, we only need to show that
$\tau_1(\mu)$ is strictly decreasing in $(0, \infty)$.

Let $\tau_1(\mu)$ be the first zero of the solution $u$ of the initial value problem
$$
\left\{
\aligned
&-(r^{N-1}u'(r))'&=\mu r^{N-1} a(r) u(r),   \ \ \ \ \qquad \ \ \   r\in (0,\infty),\\
&u'(0)=0,\\
&u(0)=\zeta.\\
\endaligned
\right.
\eqno (2.13)
$$
Let $\tau_1(\mu^*)$ be the first zero of the solution $v$ of the initial value problem
$$
\left\{
\aligned
&-(r^{N-1}v'(r))'=\mu^* r^{N-1} a(r) v(r),   \ \ \ \ \qquad \ \ \   r\in (0,\infty),\\
&v'(0)=0,\\
&v(0)=\zeta.\\
\endaligned
\right.
\eqno (2.14)
$$
We only need to show that
$$\tau_1(\mu)>\tau_1(\mu^*)\ \ \ \ \text{if}\ \mu^*>\mu.
\eqno (2.15)
$$

Suppose on the contrary that $\tau_1(\mu)\leq \tau_1(\mu^*)$. Then
$$v(r)>0, \ \  r\in [0, \tau_1(\mu)); \ \ \ \ v'(\tau_1(\mu))<0.
\eqno (2.16)$$
Multiplying the equations in (2.13) and (2.14) by $v$ and $u$, respectively, and integrating from $0$ to $\tau_1(\mu)$, we get
$$-(\tau_1(\mu))^{N-1}v(\tau_1(\mu))u'(\tau_1(\mu))
=(\mu-\mu^*)\int^{\tau_1(\mu)}_0 r^{N-1}a(r)u(r)v(r)dr.
$$
However, this is impossible from (2.16) and the fact
$$u(r)>0, \ \  r\in [0, \tau_1(\mu)); \ \ \ \ u'(\tau_1(\mu))<0.
$$
Therefore,  (2.15) is valid. \hfill{$\Box$}
 \vskip 2mm

 (2) Using the similar method to treat (2.15) and  the fact $\tau_k(\mu)$ is strictly decreasing in $(0, \infty)$, it is not difficult to show that $r_k(\mu)$ is strictly decreasing for $\mu\in (0, \infty)$.

Let $u$ be the solution of (2.13) and  $r_k(\mu)$ be the $k$-th positive zero of $u'$.
Then
$$\tau_k(\mu)<r_k(\mu)<\tau_{k+1}(\mu).$$
Without loss of generality, we may assume that
$$
u'(r)<0, \ \  r\in (r_{k-1}(\mu), r_{k}(\mu));   \ \ \ \
u(r)<0, \ \  r\in (\tau_{k}(\mu), r_{k}(\mu)).
\eqno (2.17)
$$
(The other cases can be proved by the similar method.)
Let $v$ be the solution of (2.14) and  $r_k(\mu^*)$ be the $k$-th positive zero of $v'$. Then it follows from  (2.17) that
$$v'(r)<0, \ \ \ \  r\in (r_{k-1}(\mu^*), r_{k}(\mu^*)).
\eqno (2.18)
$$

Suppose on the contrary that there exist some $k$ and some $\mu, \mu^*$ with $\mu<\mu^*$, such that
$$r_k(\mu)\leq r_k(\mu^*).
\eqno (2.19)$$
 Combining this with the fact that $\tau_k(\mu)$ is strictly decreasing in $\mu$ and using (2.17), it follows that 
$$v'(r)<0, \ \ \ r\in [\tau_k(\mu),r_k(\mu)); \ \ \ \  v(r)<0, \ \ \ r\in [\tau_k(\mu),r_k(\mu)].
$$
Multiplying the equation in (2.14) by $u$ and the equation in (2.13) by $v$ and integrating from
$\tau_k(\mu)$ to $r_k(\mu)$,
we get
$$
\aligned
 &\ \ \ \ (\tau_k(\mu))^{N-1}v(\tau_k(\mu))u'(\tau_k(\mu))+(r_k(\mu))^{N-1}v'(r_k(\mu))u(r_k(\mu))\\
&=(\mu-\mu^*)\int^{r_k(\mu)}_{\tau_k(\mu)} r^{N-1}a(r)u(r)v(r)dr.\\
\endaligned
$$
This together with the signs of $u, u', v, v'$ at
$\tau_k(\mu)$ and $r_k(\mu)$ imply that (2.19) is impossible.

\hfill{$\Box$}

\vskip3mm
 \noindent{\bf Proof of Theorem 2.2}
 Let $u(r; \zeta, \mu)$ be the unique solution of (2.6). For $k\in \mathbb{N}$. Let
  $\mu^r_k$  be such that $u'(1; \zeta, \mu^r_k)=0$ and
  $u(r; \zeta, \mu^r_k)$ has exactly $k$ zeros in $(0, 1)$.

  Let
  $$\psi_k(r):=u(r; \zeta, \mu^r_k), \ \ \ \ r\in [0, 1].
  $$
  Then Lemmas 2.1-2.5 guarantee the desired results. In particular,
   $$\psi_0(r)\equiv \zeta;\ \ \ \  \ \ \psi_1(r)\ \text{is monotone on}\ r \in (0, 1). $$
  \hfill{$\Box$}

\vskip 3mm



\noindent{\bf Lemma 2.6} Let $\{(\mu_n,y_n)\}$ be a sequence of solutions of the problem
$$-(r^{N-1}y'_n)'=\mu_n r^{N-1}g(y_n), \ \ \ \  \ y_n'(0)=y_n'(1)=0,
\eqno (2.20)
$$
where $|\mu_n|\leq \hat \mu$\,($\hat \mu$ is a positive constant), \ $g:\mathbb{R}\to \mathbb{R}$ satisfies
$$|g(s)|\leq L_0|s|\ \ \  \text{for some constant}\ L_0>0.
$$
Then $||y'_n||_\infty\to \infty$ as $n\to \infty$ implies $||y_n||_\infty\to \infty$ as $n\to \infty$.

\noindent{\bf Proof.} Assume on the contrary that $||y_n||_\infty\not\to \infty$ as $n\to \infty$.
Then, after taking a subsequence
and relabeling, if necessary, it follows that
$$||y_n||_\infty\leq M_0
\eqno (2.21)
$$
for some $M_0>0$.
From (2.20), we get
$$y_n'(r)=-\mu_n\int^r_0 \big(\frac sr\big)^{N-1}g(y_n(s))ds,
$$
which implies that
$$||y_n'||_\infty\leq \hat \mu L_0\cdot||y_n||_\infty\leq\hat \mu M_0L_0.
$$
However, this is a contradiction. \hfill{$\Box$}

 \vskip 3mm

\section{Functional setting and preliminary properties}

   The main point to prove Theorem 1.1 consists in using the unilateral global bifurcation theorem of [4, 5, 8]

Let $E$ be a real Banach space with norm $||\cdot||$. $\mathscr{E}$ will denote $E \times \mathbb{R}$.  Let the mapping $\mathcal{G}: \mathscr{E}\to E$ satisfy

{\it Assumption} $\mathfrak{A}$: if $\mathcal{G}(0,\lambda)=0$
for $\lambda\in \mathbb{R}$, $\mathcal{G}$ is completely continuous and
$$\mathcal{G}(x, \lambda)=\lambda Lx+ H(x,\lambda),$$ where $L$ is a completely continuous linear operator on $E$ and $||H(x, \lambda)||/||x||\to 0$ uniformly on bounded subsets of $\mathbb{R}$ as $||x||\to 0$.

 \vskip 2mm

 Define $\Phi(\lambda):E\to E$ by $\Phi(\lambda)(x)= x-\mathcal{G}(x,\lambda)$ and define $\mathfrak{L}$ to be the closure of $\{(x, \lambda) \in \mathscr{E}:\, x = \mathcal{G}(x, \lambda),\,  x\neq 0\}$ in $\mathscr{E}$. Then  (cp.  Rabinowitz [13]) $\mathfrak{L}\cap(\{0\}\times \mathbb{R})\subseteq \{0\}\times  r(L)$, where $r(L)$ denotes the real characteristic value of $L$.
 If $\mu\in  r(L)$, define $C_\mu$ to be the component of $\mathfrak{L}$ containing $(0,\mu)$.

Assume now that $\mu\in r(L)$ such that $\mu$ has multiplicity $1$. Suppose that
$v\in E\setminus\{0\}$ and $l\in E^*$ such that
$$v=\mu Lv, \ \ \ \ \ l=\mu L^* l,$$
(where $L^*$ is the adjoint of $L$) and $l(v)=1$. If $y\in (0,1)$, define
$$K_y=\{(u, \lambda)\in \mathscr{E}:\, |l(u)|> y||u||\},
$$
$$K^+_y=\{(u, \lambda)\in \mathscr{E}:\, l(u)> y||u||\}, \ \
  K^-_y=\{(u, \lambda)\in \mathscr{E}:\, l(u)< -y||u||\}.
$$
By [13, Lemma 1.24], there exists an $S>0$ such that
 $$(\mathcal{L}\setminus\{(\mu, 0)\}) \cap \bar {\mathscr{E}}_S(\mu)\subseteq K_y,$$
where ${\mathscr{E}}_S(\mu)=\{(u, \lambda)\in \mathscr{E}\,|\,\|u\|+|\lambda-\mu|<S\}$ and  $\bar {\mathscr{E}}_S(\mu)$ denotes closure of  $\mathscr{E}_S(\mu)$.
 For $0<\epsilon\leq S$ and $\nu=\pm$, define $D^\nu_{\mu, \epsilon}$ to be the component of
 $\{(0, \mu)\}\cup (\mathfrak{L}\cap\bar {\mathscr{E}}_\epsilon(\mu)\cap K_y^\nu)$
 containing $(0, \mu)$, $C^\nu_{\mu, \epsilon}$ to be the component of
 $\overline{C_\mu \setminus \mathcal{D}^{-\nu}_{\mu, \epsilon}}$ containing $(0, \mu)$
 (where $-\nu$ is interpreted in the natural way), and $C_{\mu,\nu}$ to be the closure of
 $\bigcup_{S\geq \epsilon >0} C^\nu_{\mu, \epsilon}$. Then $C_{\mu,\nu}$ is connected and,
 by [5], $C_\mu=C_{\mu,+}\cup C_{\mu,-}$. By [13, Lemma 1.24],  the definition of $C_{\mu,\nu}$ is independent of $y$.
 \vskip 3mm

 \noindent{\bf Theorem 3.1} [5, Theorem 2] Either $C_{\mu,+}$ and $C_{\mu,-}$ are both unbounded or
 $$C_{\mu,+} \cap C_{\mu,-}\neq\{(0, \mu)\}.$$

\vskip 3mm

\section{Proof of the Main Results}

Let $X:=\{u\in C^1[0,1]\,|\, u'(0)=u'(1)=0\}$. Then it is a Banach space under the norm
$$||u||_X=\max\{||u||_\infty, \ ||u'||_\infty\}.$$

\vskip 3mm

We shall prove that the first choice of the alternative of Theorem 3.1 is the only possibility.

 In what follows, we use the terminology of Rabinowitz [13]. Let $S_{k,+}$ denote the set of functions in $X$
which have exactly $k-1$ interior nodal (i.e. non-degenerate) zeros in $(0, 1)$ and are positive near
$r = 0$, set $S_{k,-} =-S_{k,+}$ , and $S_k = S_{k,+} \cup  S_{k,-}$.  Finally, let $\Phi_{k, \pm} = \mathbb{R}\times S_{k, \pm}$ and $\Phi_k = \mathbb{R}\times S_k$  under the product topology.

\vskip 3mm

Let us consider the problem
$$
\left\{
\aligned
-\Delta u+u= f(u)    \ \ \ \ \qquad & \text{in}\ B_1,\\
u>0\qquad\qquad\quad\quad \ \ \ \ \           & \text{in}\ B_1,\\
\partial_\nu u=0 \ \ \ \qquad\quad\qquad \ \ \ &\text{on}\ \partial B_1, \\
\endaligned
\right.
\eqno (4.1)
$$
which is equivalent to
$$
\left\{
\aligned
&-u''-\frac {N-1}r u'+u= f(u), \ \ \ \ r\in (0,1),\\
&u>0,\ \ \ \ \ \ \ \ \ \ \  \ \ \ \ \  \ \ \quad \ \qquad \qquad r\in [0, 1], \\
& u'(0)=u'(1)=0.\\
\endaligned
\right.
\eqno (4.2)
$$
Let
$$v:=u-\beta.
$$
Then (4.2) can be rewritten as
$$
\left\{
\aligned
&-v''-\frac {N-1}r v'+v=f(v+\beta)-\beta, \ \ \ \ r\in (0,1),\\
&v>-\beta,\ \ \ \ \ \ \ \ \ \ \  \ \ \   \ \quad \qquad \ \quad \ \qquad \qquad r\in [0, 1], \\
& v'(0)=v'(1)=0.\\
\endaligned
\right.
\eqno (4.3)
$$
Let
$$h(s):=
\left\{
\aligned
f(s+\beta)-\beta,\ \ \ \  \ &s\geq -\beta,\\
-\beta, \ \qquad \qquad\  \ \  \ &s<-\beta.\\
\endaligned
\right.
$$
Then
$$h(v)=h'(0)v+\xi(v)=f'(\beta)v+\xi(v),\ \ \ \ h'(0)=f'(\beta),
$$
and
$$\xi'(0):=\lim_{v\to 0}\frac{\xi(v)}{v}=0.
\eqno (4.4)
$$
Thus, to study the $S_{k,\nu}$-solutions of (4.3), let us consider the auxiliary problem
$$
\left\{
\aligned
&-v''-\frac {N-1}r v'+v=\lambda f'(\beta)v+\lambda\xi (v), \ \ \ \ r\in (0,1),\\
&v>-\beta,\ \ \qquad \ \  \ \ \ \ \ \ \ \ \   \ \quad \qquad \ \quad \ \qquad \qquad r\in [0, 1],  \\
& v'(0)=v'(1)=0.
\endaligned
\right.
\eqno (4.5)
$$

For $e\in X$, let $Te$ be the unique solution of the problem
$$
\left\{
\aligned
&-z''-\frac {N-1}r z'+z=e, \ \ \ \  r\in (0,1),\\
& z'(0)=z'(1)=0.
\endaligned
\right.
$$
Then the map $T:X\to X$ is completely continuous, and
$$r(T)=\big\{\lambda^r_j\,|\, \lambda^r_j=\mu^r_{j-1}+1, \; j=1,2, \cdots\big\}.
$$
Here $r(T)$ denotes the real characteristic value of $T$.
Obviously£¬ (4.5) is equivalent to
$$
v=\lambda T(f'(\beta)v)+\lambda T\xi(v),
\eqno (4.7)
$$
$$
v>-\beta.\hskip 27mm
\eqno (4.8)
$$

To show that (4.5) has a
 $S_{k,\nu}$-solution, let us consider the auxiliary problem (4.7)-(4.8)
as a bifurcation problem from the trivial solution $v\equiv 0$.
Furthermore, we have from (4.4)  that
$$
\frac{||T \xi(v)||_{X} }{||v||_{X}}
\leq ||T||_{X\to X}
               \max\Big\{\frac{||\xi(v)||_\infty}{||v||_X},\frac{||\xi'(v)v'||_\infty }{||v||_X}\Big\}
\to 0\ \ \  \text{as}\ ||v||_{X}\to 0.
$$
\vskip 3mm

Now the Dancer's unilateral global bifurcation theorem
for (4.7) can be stated as follows:

     Let
     $$\mathfrak{L}:=\overline{\big\{(\lambda, v)\in (0, \infty)\times X: \, (\lambda, v) \, \text{satisfies}\ (4.7), \, v\neq 0\big\}}^X.$$
     For $\lambda^r_k \in  r(T)$, define $C_k$ to be the component of $\mathfrak{L}$ containing $\big(\frac{\lambda^r_k}{f'(\beta)},0\big)$. Then
      $$C_k:= C_{k,+} \cup C_{k,-},
           $$
      where
      $$C_{k,\nu}:= C_{{\lambda{{_k^r}/f'(\beta)}},\,\nu} \ \  \ \ \ \nu\in \{+, -\},$$
      see Section 3 for detail. Now the Dancer's unilateral global bifurcation theorem yields that  either $C_{k,+}$ and $C_{k,-}$ are both unbounded or
 $$C_{k,+} \cap C_{k,-}\neq\{\big(\frac{\lambda^r_k}{f'(\beta)},0\big)\}.
 \eqno (4.9)
 $$

\vskip 3mm

From (A1), it follows that if $(\lambda, v)$ is a solution of
$$
\left\{
\aligned
&-v''-\frac{N-1}r v'+v=\lambda h(v),\ \\
&v(\tau)=v'(\tau)=0\\
\endaligned
\right.
$$
for some $\tau\in (0, \infty)$, then $v\equiv 0$. This implies that
$$C_{k,+} \subset \Big(\Phi_{k,+}\cup \big\{\big(\frac{\lambda_k^r}{f'(\beta)}, 0\big)\big\}\Big).
$$

Clearly, if (4.9) holds, then there exists $(\lambda_\ast, v_\ast)\in C_{k,+}\cap C_{k,-}$, such that 
$(\lambda_\ast, v_\ast)\neq\big(\frac{\lambda^r_k}{f'(\beta)},0\big),$  and $v_\ast\in S_{k,+}\cap S_{k,-}$, which contradicts the definition of $S_{k,+}$ and $S_{k,-}$.

Furthermore, we get
\vskip 3mm

\noindent{\bf Lemma 4.1.}\ For given $k\geq 2$, $C_{k,+}$ and $C_{k,-}$ are both unbounded, and  $\big(\frac{\lambda_k^r}{f'(\beta)}, 0\big)$  bifurcates two unbounded components $C_{k,+}$ and $C_{k,-}$ of solutions to problem (4.7), such that
$$\big(C_{k,+}\setminus \big\{\big(\frac{\lambda_k^r}{f'(\beta)}, 0\big)\big\}\big)\subseteq \Phi_{k,+}, \ \ \ \ \  \ \big(C_{k,-}\setminus \big\{\big(\frac{\lambda_k^r}{f'(\beta)}, 0\big)\big\}\big)\subseteq \Phi_{k,-}.
$$

\vskip 3mm
\noindent{\bf Lemma 4.2.}\  Let $(\lambda, v)\in C_{k,\nu}$ with $\lambda\in [0,1]$. Then
$$v(r)>-\beta, \ \ \ \ r\in [0, 1].
\eqno (4.10)
$$

\noindent{\bf Proof.} Suppose on the contrary that there exists $x_0\in[0,1]$ such that
$$v(x_0)=\underset{r\in [0,1]}\min\;  v(r)=-\beta.
$$
Then there exists $r_0\in [0,1]$ such that either
$$ v(r_0)=0, \ \  v(r)<0\ \text{for} \ r\in [x_0, r_0), \ \  v'(r)>0\ \text{for} \ r\in (x_0, r_0];
\eqno (4.11)
$$
or
$$ v(r_0)=0, \ \  v(r)<0\ \text{for} \ r\in (r_0,x_0], \ \  v'(r)<0\ \text{for} \ r\in [r_0, x_0).
\eqno (4.12)
$$

   We only deal with the case (4.11), the case (4.12) can be treated by the similar way.

    By (A1)-(A3), there exists $m\geq 0$ such that $h(s)+ ms$ is monotone increasing in $s$ for $s\in [-\beta, +\infty)$. Then
$$-v''-\frac {N-1}r v'+v+\lambda mv=\lambda[h(v)+mv], \ \ \   \ r\in (0, 1],
$$
and, since
$$- (-\beta)''-\frac {N-1}r(-\beta)' +(-\beta)+\lambda m(-\beta)\leq\lambda[h(-\beta)+m(-\beta)],\ \ \   \ r\in (0,1],$$
it follows that
$$-(v+\beta)''-\frac {N-1}r(v+\beta)' +(\lambda m+1)(v+\beta)\geq \lambda\big([h(v)+mv]-[h(-\beta)+m(-\beta)]\big)\geq 0, \
\ \ \ \ \   \ r\in (0, 1].
$$
Denote
$$w:=v+\beta.
$$
Then
$$
\aligned
&w''+\frac {N-1}r w'-(\lambda m+1)w \leq 0,\ \ \quad \ \
\    \ r\in (0, 1],\\
&w'(0)=w'(1)=0.\\
\endaligned
$$
It follows from [6, Theorem 3.5] or [12, Theorem 3 in Chapter 1] that, $w$ cannot achieve a non-positive minimum in the interval $(0, 1)$ unless it is constant. From (4.11), it follows that
$$\inf_{[x_0, r_0]} w(r)=\min\{w(x_0), w(r_0)\}=w(x_0)=0.
$$
This together with $w'(x_0)=0$ imply that
$$w(r)\equiv 0, \ \ \ \ \ r\in [x_0, r_0].$$
However, this contradicts the fact $w'(r)>0,\ r\in(x_0, r_0)$.
Therefore,
$$v(r)>-\beta, \ \ \ \ \ r\in [0, 1].$$
\hfill{$\Box$}

\vskip 4mm

In view of  Lemma 4.2, (4.5) is equivalent to (4.7). So, we only need to show that
$$C_{k,\nu}\cap\big(\{1\}\times X\big)\neq \emptyset.
\eqno (4.13)
$$
In the following, we only deal with the case ` $\nu=+$ ' since the other case can be treated by the similar way.

\vskip 3mm

 Let $k\geq 2$ be fixed, and let
$(\eta_n, y_n)\in C_{k,+}$ satisfy
$$\eta_n+||y_n||_{X}\rightarrow
\infty.
$$
It is easy to check that
$$\eta_n>0, \ \ \ \ n\in \mathbb{N}.
\eqno (4.14)$$
From (A3), it follows that that \ $\frac{\lambda_k^r}{h'(0)}<1$, i.e.
$$\frac{\lambda_k^r}{f'(\beta)}<1.
\eqno (4.15)
$$

We shall show that
$$ C_{k,+}\cap\big(\{1\}\times X\big)\neq \emptyset.
\eqno (4.16)
$$

Assume on the contrary that $ C_{k,+}\cap\big(\{1\}\times X\big)= \emptyset$. Then
 $$C_{k,+}\subset (0, 1)\times X,
 $$
 and accordingly,
 $$0<\eta_n<1.
 $$
 Thus
$$||y_n||_X\to \infty,\ \ \ \ \ \ n\to\infty,
\eqno (4.17)
$$
which together with Lemma 2.6 imply that
$$||y_n||_\infty\to \infty,\ \ \ \ \ \ n\to\infty.
\eqno (4.18)
$$
This means that $C_{k,+}$ is unbounded in $C[0,1]$!

We may assume that $\eta_n\to \bar\eta\in [0, 1]$ as $n\to\infty$.
Let
$$z_n:=\frac {y_n}{||y_n||_\infty}.
$$
Then  $||z_n||_\infty=1$ and
$$\left\{
\aligned
&-z_n''-\frac {N-1}r z_n'+z_n=\eta_n \frac{h(y_n)}{y_n}z_n, \ \ \ \ r\in (0,1),\\
& z_n'(0)=z_n'(1)=0.\\
\endaligned
\right.
\eqno (4.19)
$$
From (A1)-(A3) and the definition of $h$, it follows that $\frac{h(y_n(r))}{y_n(r)}$ is continuous in $[0,1]$ and is bounded  uniformly in $n$. After taking subsequence if necessary, we may assume that
$$(\eta_n, z_n)\to (\bar \eta, z^*), \ \ \ \  \text{in} \ \mathbb{R}\times X.
\eqno (4.20)
$$
 Here $||z^*||_\infty=1$.

 As a direct consequence of the Banach contraction mapping principle in a small neighborhood of $\tau$,  the initial value problem
$$\left\{
\aligned
&-{z}''-\frac {N-1}r {z}'+z=\bar\eta H(r)z,\\
& z(\tau)={z}'(\tau)=0\\
\endaligned
\right.
$$
has a unique solution $z\equiv 0$. Notice that taking subsequence if necessary, we may assume that $\frac{h(y_n)}{y_n}\overset{w}\longrightarrow H$ in $L^2[0,1]$.
So, all of zeroes of $z^*$ are simple,
and accordingly $(\bar \eta, z^*)\in C_{j,+}$ for some $j\in \mathbb{N}$.
\vskip 3mm

Let
$$\tau(1,n)<\cdots<\tau(k-1,n)
$$
denote the zeros of $y_n$, and let
$$\tau(0,n):=0, \ \ \ \ \tau(k,n):=1.
$$
 Then, after taking a subsequence if
necessary,
$$
\lim_{n\rightarrow \infty}\tau(l,n):=\tau(l,\infty), \qquad l\in
\{0, 1, \cdots, k-1, k \}. \eqno (4.21)
$$
Denote
$$J_l:=\big(\tau(l,\infty), \tau(l+1,\infty)\big), \ \ \ \ l\in
\{0, 1, \cdots, k-1 \}.
$$

 {\bf Claim }\  We claim that
$$J_l=\emptyset\ \ \  \qquad \text{if}\ l\in
\{0, 1, \cdots, k-1 \} \ \text{and} \ l \text{\ is odd},
\eqno (4.22) $$
and
$$\underset{n\to\infty}\lim y_n(r)=+\infty\ \text{uniformly in}\ [\tau(l,\infty)+\epsilon, \tau(l+1,\infty)-\epsilon] \ \  \text{if}\ l\in
\{0, 1, \cdots, k-1 \} \ \text{and} \ l \text{\ is even},
\eqno (4.23)
$$
where $\epsilon>0$ is small constant.

In fact,
suppose on the contrary that
$$J_{l_0}\neq\emptyset\ \ \  \qquad \text{for some}\ l_0\in
\{0, 1, \cdots, k-1 \} \ \text{and} \ l_0 \text{\ is odd}.
$$
Then  we have from Lemma 4.2 that
$$-\beta <y_n(r)< 0,\quad \ \quad  \; r\in (\tau(l_0,n),\tau(l_0+1,n)).
$$
Thus, for any $r\in (\tau(l_0,n),\tau(l_0+1,n))$, it follows from (4.18)
and
$$-z_n''-\frac {N-1}r z_n'+z_n=\eta_n \frac{h(y_n)}{||y_n||_\infty}, \ \ \ \ r\in (\tau(l_0,n),\tau(l_0+1,n)),
$$
that
$$\left\{
\aligned
&-{z^*}''-\frac {N-1}r {z^*}'+z^*=0, \ \ \ \ r\in J_{l_0},\\
& {z^*}(\tau)={z^*}'(\tau)=0,\\
\endaligned
\right.
$$
for some $\tau\in J_{l_0}$. This implies that
$$z^*(r)= 0 = {z^*}'(r),\ \ \ \ r\in J_{l_0}.
$$
However, this contradicts the fact the solution
$(\bar \eta, z^*)\in C_{j,+}\subset [0, 1]\times S_{j,+}$ for some $j\in \mathbb{N}$.
Therefore, (4.22) is true!

Obviously, (4.23) is an immediate consequence of the fact that all of the zeros of $z^*\in S_{j,+}$ are simple and $l\in\{0,1, \cdots, k-1\}$ is even.

Therefore, the {\bf Claim} is true!

\vskip 3mm

In the following, we shall use some idea from the proof of [7, Lemma 3.2] and the proof of main results of [3, 9] to show (4.16) is valid.

\vskip 3mm

Let  $(y_n)^-$ be  the negative part of $y_n$.
Then it follows from Lemma 4.2 that
 $0\leq (y_n)^-< \beta$ since $\eta_n\in (0,1)$,
 and consequently,
$$(z_n)^-\, \to \, 0, \ \ \ \ \text{as}\ n\to \infty.
$$
Combining this with the {\bf Claim} and using the definition of $h$, it concludes that
$$\left\{
\aligned
&-{z^*}''-\frac {N-1}r {z^*}'+z^*=\bar \eta f_{+\infty} {(z^*)}^+, \ \ \ \text{a.e.} \ r\in (0,1),\\
& {z^*}'(0)={z^*}'(1)=0,\\
\endaligned
\right.
\eqno (4.24)
$$
where ${(z^*)}^+$ is the positive part of $z^*$. Now,
it follows from [6, Theorem 3.5] or [12, Theorem 3 in Chapter 1] that, $z^*$ cannot achieve a non-positive minimum in the interval $(0, 1)$ unless it is constant. Since
$z_n(0)>0$, we get
$$z^*(0)\geq 0.
$$
If $z^*(0)=0$, then it follows from (4.24) that
 $$z^*\equiv 0,\ \ \ \ \ \ \ \ r\in (0,1).$$
 However, this contradicts $||z^*||_\infty= 1$.
  So
$$z^*(0)>0.
$$
Again, from [6, Theorem 3.5] or [12, Theorem 3 in Chapter 1] that $z^*(r)>0$ in $[0,1]$. This means that $z^*\in S_{1,+}$, and therefore, since $S_{1,+}$ is open and $||z_n-z^*||_X\to 0$ as $n\to \infty$, $z_n\in S_{1,+}$ for $n$ large enough. However, this contradicts
$z_n\in S_{k,+}$ for all $n\in \mathbb{N}$ and $k\geq 2$.

Therefore, (4.16) is valid, and we may take $v_{k,+}\in \big(C_{k,+}\cap\big(\{1\}\times X\big)\big)$. Similarly, we may take $v_{k,-}\in \big(C_{k,-}\cap\big(\{1\}\times X\big)\big)$.

\vskip 3mm

To show that $v_{2,+}$ is decreasing in $[0,1]$. Let us denote $t_1\,(0<t_1<1)$
be the zero of $v_{2,+}$. Notice that $v_{2,+}$ satisfies (4.3),  i.e.,
$$
\left\{
\aligned
&-(r^{N-1}(v{_{2,+}})')'=r^{N-1}[f(v{_{2,+}}+\beta)-(v{_{2,+}}+\beta)], \ \ \ \ r\in (0,1),\\
&v{_{2,+}}>-\beta,\\
& v{_{2,+}}'(0)=v{_{2,+}}'(1)=0.\\
\endaligned
\right.
\eqno (4.25)
$$
Combining this with (A4)  and using Lemma 4.2, it concludes that
$$(r^{N-1}(v{_{2,+}})')'<0, \ \ \ \ t\in (0, t_1); \ \ \ \ \
(r^{N-1}(v{_{2,+}})')'>0, \ \ \ \ t\in (t_1,1).
\eqno (4.26)
$$
This together with the boundary condition $ v{_{2,+}}'(0)=v{_{2,+}}'(1)=0$ imply that
$$(v{_{2,+}})'<0, \ \ \ \ t\in (0, t_1); \ \ \ \ \
(v{_{2,+}})'<0, \ \ \ \ t\in (t_1,1).
$$
Therefore, $v{_{2,+}}$ is decreasing in $[0,1]$.

Using the same method, with the obvious changes, we may deduce that  $v{_{2,-}}$ is increasing in $[0,1]$.
\hfill{$\Box$}

\vskip 1cm

\centerline {\bf REFERENCES}\vskip5mm\baselineskip 20pt
\begin{description}

\item{[1]}  D. Bonheure,   B. Noris, T.  Weth, Increasing radial solutions for Neumann problems without growth restrictions, Ann. Inst. H. Poincar\'{e} Anal. Non Lin\'{e}aire 29(4) (2012), 573-588.

\item{[2]}  D.  Bonheure, E. Serra, P. Tilli, Radial positive solutions of elliptic systems with Neumann boundary conditions, J. Funct. Anal. 265(3) (2013), 375-398.

\item{[3]} G. Dai, R. Ma, Unilateral global bifurcation phenomena and nodal solutions for $p$-Laplacian, J. Differential Equations 252(3) (2012), 2448-2468.

\item{[4]} E. N. Dancer, Bifurcation from simple eigenvalues and eigenvalues of geometric multiplicity one, Bull. Lond. Math. Soc. 34 (2002), 533-538.

\item{[5]} E. N. Dancer, On the structure of solutions of non-linear eigenvalue problems,
Indiana Univ. Math. J. 23 (1973/74), 1069-1076.

\item{[6]} D. Gilbarg, \  N. S. Trudinger,  \ Elliptic Partial Differential Equations of Second Order, Springer-Verlag, Berlin, 2001, Reprint of the 1998 edition.

\item{[7]}  A. C. Lazer, P. J. Mckenna, Global bifurcation and a theorem of Tarantello,
J. Math. Anal. Appl. 181 (1994), 648-655.

\item{[8]} J. L\'{o}pez-G\'{o}mez, Spectral Theory and Nonlinear Functional Analysis, in: Research Notes in Mathematics, vol. 426, Chapman \& Hall/CRC, Boca Raton, Florida, 2001.

\item{[9]}  R. Ma, B. Thompson, Nodal solutions for nonlinear eigenvalue problems, Nonlinear Anal. 59 (2004), 707-718.

\item{[10]}\ A. Miciano, R. Shivaji, Multiple positive solutions for
a class of semipositone Neumann two point boundary value problems,
J. Math. Anal. Appl. 178 (1993), 102-115.

\item{[11]}   L. A. Peletier, J. Serrin, Uniqueness of positive solutions of semilinear equations in $\mathbb{R}^n$ . Arch. Rational Mech. Anal. 81(2) (1983),  181-197.

\item{[12]}   M. H. Protter, H. F. Weinberger, Maximum Principles in Differential Equations, Springer-Verlag, New York, 1984.

\item{[13]} P. H. Rabinowitz, Some global results for nonlinear eigenvalue problems, J. Funct. Anal. 7 (1971), 487-513.

\item{[14]}  E. Serra, P. Tilli, Monotonicity constraints and supercritical Neumann problems, Ann. Inst. H. Poincar\'{e} Anal. Non Lin\'{e}aire 28 (2011), 63-74.

\item{[15]} W. Walter, Ordinary Differential Equations. Translated from the sixth German (1996) edition by Russell Thompson. Graduate Texts in Mathematics 182.  Springer-Verlag, New York, 1998.

\end{description}
\end{document}